\documentclass[secthm,seceqn,amsthm,ussrhead,12pt]{amsart}
\usepackage{amsmath,latexsym}
\usepackage[english]{babel}
\usepackage[psamsfonts]{amssymb}
\usepackage{times}
\usepackage{cite}

\usepackage[mathcal]{euscript}
\numberwithin{equation}{section} \textwidth=17.5cm
\newtheorem{theorem}{Theorem}[section]
\newtheorem{lemma}[theorem]{Lemma}

\newtheorem{corollary}[theorem]{Corollary}
\newtheorem{defin}[theorem]{Definition}

\numberwithin{equation}{section}

\begin{document}

\title{Geometric Description of the Preduals
of Atomic Commutative von Neumann Algebras}

 \author[M. M. Ibragimov]{M. M. Ibragimov}
 \author[K. K. Kudaybergenov]{K. K. Kudaybergenov}
 \author[S. J.  Tleumuratov]{S. J. Tleumuratov}
\author[J. X.  Seypullaev]{J. X. Seypullaev}

\address{Department of Mathematics,
 Karakalpak state university
Ch. Abdirov 1
  230113, Nukus,    Uzbekistan,}

\email{mukhtar$_{-}$nukus@mail.ru}

\email{karim2006@mail.ru}

\email{sarsenmat@rambler.ru}

\email{jumabek81@mail.ru}
\begin{abstract}
Strongly facially symmetric spaces isometrically isomorphic to the
predual space of an atomic commutative von Neumann algebra are
described.
\end{abstract}


\maketitle
\section{Introduction}

An important problem of the theory of operator algebras is a
geometric characterization of state spaces of operator algebras.
In the mid-1980s, Friedman and Russo wrote the paper~\cite{fr1}
related to this problem, in which they introduced facially
symmetric spaces, largely for the purpose of obtaining a geometric
characterization of the predual spaces of $JBW^\ast$-triples
admitting an algebraic structure. Many of the properties required
in these characterizations are natural assumptions for state
spaces of physical systems. Such spaces are regarded as a
geometric model for states of quantum mechanics. In \cite{fr3}, it
was proved that the preduals of von Neumann algebras and, more
generally, $JBW^\ast$-triples are neutral strongly facially
symmetric spaces.

The project of classifying facially symmetric spaces was initiated
in \cite{fr4}, where a geometric characterization of complex
Hilbert spaces and complex spin factors was given. The
$JBW^\ast$-triples of ranks $1$ and $2$ and Cartan factors of
types $1$ and $4$ were also described. Afterwards, Friedman and
Russo obtained a description of atomic facially symmetric spaces
\cite{fr5}. Namely, they showed that a neutral strongly facially
symmetric space is linearly isometric to the predual of one of the
Cartan factors of types 1--6 provided that it satisfies four
natural physical axioms, which hold for the predual spaces of
$JBW^\ast$-triples. In the 2004 paper \cite{NR}, Neal and Russo
found geometric conditions under which a facially symmetric space
is isometric to the predual of a $JBW^\ast$-triple. In particular,
they proved that any neutral strongly facially symmetric space
decomposes into a direct sum of atomic and nonatomic strongly
facially symmetric spaces. This paper describes strongly facially
symmetric spaces isometrically isomorphic to preduals of atomic
commutative von Neumann algebras.

This paper describes strongly facially symmetric spaces
isometrically isomorphic to preduals of atomic commutative von
Neumann algebras.

\section{Preliminaries}

  Let  $Z$ be a real or complex normed space. We say that elements
    $x, y \in Z$
  are \textit{orthogonal} and write  $x \diamondsuit y$ if
$$
\|x + y\|=\|x - y\|=\|x\|+\|y\|.
$$
We say that subsets   $S, T\subset Z$  are \textit{orthogonal} and
write
  $(S \diamondsuit T)$ if  $x
\diamondsuit y$ for all  $(x, y)\in S\times T.$ For a subset $S$
of $Z,$ we put
$$
S^\diamondsuit =\{x \in Z: x \diamondsuit y,\,
\forall\,  y\in S\};
$$ the set  $S^\diamondsuit$ is called the
\textit{orthogonal complement} of $S.$ A convex subset $F$  of the
unit ball
$$
Z_1=\{x\in Z:\|x\|\leq 1\}
$$ is called a \textit{face} if the
relation
 \begin{center}
 $
  \lambda y+(1-\lambda)z\in F,
 $
 where
 $y, z\in Z_1,$ $\lambda
\in (0,1),$
\end{center}
 implies   $y, z\in F.$ A face
 $F$ of the unit ball is said to be \textit{norm exposed} if
   $$
   F=F_u=\{x\in Z: u(x)=1\}
   $$
   for some  $u\in Z^*$ with  $\|u\| =1.$
An element   $u\in Z^\ast$ is called a projective unit if  $\|u\|
=1$ and  $u(y)=0$ for all   $y\in F_u^\diamondsuit$ (see \cite
{fr1}).

\begin{defin}  (\cite{fr1}). A norm exposed face $F_{u}$ in $Z_1$
is called a \textit{symmetric face} if there exists a linear
isometry $S_u$ from $Z$ to $Z$ such that  $S_u^2=I$ whose fixed
point set coincides with the topological direct sum of the closure
$\overline{sp}F_{u}$ of the linear hull of the face $F_{u}$ and
its orthogonal complement $F_{u}^{\diamond},$ i.e., with
$(\overline{sp}F_{u})\oplus F_{u}^\diamond.$
\end{defin}

\begin{defin} (\cite{fr1}). A space $Z$ is said to be \textit{weakly facially symmetric}
  $(WFS)$ if each norm exposed
face in  $Z_1$ is symmetric.
\end{defin}

For each symmetric face  $F_u,$ contractive projections
$P_k(F_u),$ $k = 0, 1, 2,$  on  $Z$ are defined as follows. First,
$P_1(F_u) = (I - S_u)/2$  is the projection onto the eigenspace
corresponding to the eigenvalue $-1$ of the symmetry  $S_u.$ Next,
 $P_2(F_u)$ and  $P_0(F_u)$ are defined as projections of  $Z$ onto
$\overline{sp}F_u$ and   $F_u^\diamondsuit,$ respectively; i.e.,
$P_2(F_u) + P_0(F_u) = (I + S_u )/2.$ The projections  $P_k(F_u)$
are called the \textit{geometric Peirce projections}.

A projective unit  $u$  from  $Z^*$  is called a \textit{geometric
tripotent} if  $F_u$ is a symmetric face and $S_u^*u=u$  for the
symmetry $S_u$  corresponding to
 $F_u.$  By  $\mathcal{GT}$ and
$\mathcal{SF}$  we denote the sets of all geometric tripotents and
symmetric faces, respectively; the correspondence
 $$
 \mathcal{GT}
\ni u \mapsto F_u \in \mathcal{SF}
$$  is one-to-one (see
\cite[Proposition 1.6]{fr2}).  For each geometric tripotent $u$
from the dual $WFS$ space $Z,$  we denote the Peirce projections
by
$$
P_k(u) = P_k(F_u), k = 0, 1, 2.
$$
 We set
 $$
 U= Z^\ast, Z_k(u) =
Z_k(F_u)= P_k(u)Z,\, U_k(u) = U_k(F_u)= P_k(u)^\ast(U).
$$
The Peirce decomposition
  $$
  Z = Z_2(u) + Z_1(u) + Z_0(u),\,
  U =
U_2(u) + U_1(u) + U_0(u)
$$
holds. Tripotents  $u$ and   $v$   are said to be
\textit{orthogonal} if $u \in U_0(v)$ (which implies  $v \in
U_0(u)$)  or, equivalently,
 $u \pm v \in
\mathcal{GT}$ (see  \cite[Lemma 2.5]{fr1}). More generally,
elements $a$ and   $b$ of $U$ are said to be \textit{orthogonal}
if one of them belongs to $U_2(u)$ and the other belongs to
$U_0(u)$  for some geometric tripotent  $u.$

 A contractive projection $Q$ on  $Z$ is said to be \textit{neutral}
  if $\|Q x
\|= \|x\|$ implies   $Q x=x$ for each $x \in Z.$   A space $Z$ is
said to be \textit{neutral} if, for each symmetric face $F_u,$ the
projection $P_2(u),$  corresponding to the symmetry $S_u,$ is
neutral. A space  $Z$  is said to be \textit{atomic} if each
symmetric face of  $Z_1$ contains an extreme point.

\begin{defin} (\cite{fr1}).  A $WFS$ space $Z$
is said to be \textit{strongly facially symmetric} $(SFS)$ if, for
each norm exposed face $F_u$ of  $Z_1$ and each $g \in Z^\ast$
satisfying the conditions $\|g\| = 1$ and  $F_u \subset F_g,$ we
have $S_u^\ast g = g,$ where $S_u$ is the symmetry corresponding
to $F_u.$
\end{defin}

Instructive examples of neutral strongly facially symmetric spaces
are Hilbert spaces, the preduals of von Neumann algebras or
$JBW^\ast$-algebras, and, more generally, the preduals of
$JBW^\ast$-triples. Moreover, geometric tripotents correspond to
nonzero partial isometries of von Neumann algebras and tripotents
in $JBW^\ast$-triples (see \cite{fr3}).

 In a neutral
strongly facially symmetric space $Z,$ each nonzero element admits
a polar decomposition \cite[Theorem  4.3]{fr2},  that is, for  $0
\neq x \in Z$ there exists a unique geometric tripotent  $v =v_x$
for which $v(x) = \|x\|$ and  $\langle v, x^\diamond \rangle = 0.$
If  $x, y \in Z,$ then  $x \diamond y$ if and only if $v_x
\diamond v_y$ (see \cite[Corollary 1.3(b) and Lemma  2.1]{fr1}). A
geometric tripotent $u\in U$ is said to be minimal if $\dim
U_2(v)=1.$ Note that if $u$ is a minimal geometric tripotent, then
there exists a $z\in Z$ such that $P_2(u)(x)=u(x) z$ for all $x\in
Z.$

 A neutral $SFS$-space Z is said to
be \textit{point exposed (PE)} if each extreme point of the ball
$Z_1$  is a norm exposed point \cite{fr4}.

The set of geometric tripotents is ordered as follows: given  $u,
v \in \mathcal{GT},$ we set $u \leq v$ if  $F_u \subset F_v.$
 Note
that this is equivalent to the relation $P_2(u)^*v = u$ or to the
condition that  $v - u$ either vanishes or is a geometric
tripotent orthogonal to $u$ (see \cite[Lemma 4.2]{fr2}).

\section{Main result} \label{subsec3}

In what follows,  $Z$ is always assumed to be a neutral strongly
facially symmetric space for which there exists a geometric
tripotent $e$ satisfying the condition $P_2(e)=I,$ i.e., $Z_2
(e)=Z.$ We set
\[ \nabla = \{u \in \mathcal{GT}: u \le e\} \cup
\{0\}.
\]
 As is known \cite[Proposition 4.5]{fr2},
the set $\nabla$  is a complete orthomodular lattice with
orthocomplement
 $u^\perp = e - u$ with respect to the order
 $"\le".$

 First, we prove several auxiliary lemmas.

\begin{lemma}
\label{lem1} If $\nabla$ is a Boolean algebra, then, for any $u\in
\nabla,$ $u\neq 0,$

(a)   $P_1(u)=0;$

(b)  $P_2 (u) = P_0 (u^\perp).$
\end{lemma}

Proof. (a)  Take $u \in \nabla, u\neq 0.$ First, let us show that
\[
F_e\subseteq \overline{sp} F_u\oplus \overline{sp}F_{u^\perp}.
\]
Suppose that, on the contrary, there exists an element
\[
a\in F_e\setminus  \overline{sp} F_u\oplus \overline{sp}F_{u^\perp}.
\]
On the subspace
 $sp \{a\}\oplus  \overline{sp} F_u\oplus \overline{sp}F_{u^\perp}$ we define a functional
 $g$ by
 \[
g(\lambda a+x)=\lambda,\,\, \lambda\in \mathbb{C},\,
x\in  \overline{sp} F_u\oplus \overline{sp}F_{u^\perp}.
 \]
By the Hahn --- Banach theorem, $g$  admits an extension to $Z,$
with the same norm; we denote this extension by the same symbol
$g.$  We have
\[
\|g\|=1,\, g(a)=1,
\]
\begin{equation}
\label{fun1}
g|_{\overline{sp} F_u\oplus \overline{sp}F_{u^\perp}}=0.
\end{equation}
Take a tripotent  $v\in Z^{\ast}$  for which
 $F_g=F_v.$ By virtue of \eqref{fun1}, we have
 \[
v(a)=1,
\]
\begin{equation}\label{trip}
v|_{\overline{sp} F_u\oplus \overline{sp}F_{u^\perp}}=0.
\end{equation}

 Let us show
that $u\wedge v=0$ and $u^{\perp}\wedge v=0.$  Indeed, if $u\wedge
v \neq 0,$  then there exists an $x\in Z_1$  for which
$u(x)=v(x)=1,$  so that $x\in F_u$ and  \eqref{trip} implies that
$v(x)=0,$  which is false. Thus,
\[
u\wedge v=0, u^{\perp}\wedge v=0.
\]
 It follows from $a\in F_e\cap F_v$
 that $e\wedge v\neq 0.$  On the other hand, we have
\[
e\wedge v=(u\vee u^{\perp})\wedge v= (u\wedge v)\vee
(u^{\perp}\wedge v)=0\vee 0=0,
\]
because $\nabla$ is a Boolean algebra. This contradiction implies
\[
F_e\subseteq \overline{sp} F_u\oplus \overline{sp}F_{u^\perp}.
\]
 Since $\overline{sp}F_e=Z,$  it follows that
\begin{equation}\label{deco}
Z=\overline{sp} F_u\oplus \overline{sp}F_{u^\perp}.
\end{equation}

Now, let us show that $P_1(u)=0.$  Relation  \eqref{deco} means
that
\[
I=P_2(u)+P_2(u^{\perp}).
\]
Hence
\[
P_1(u)+P_0(u)=I-P_2(u)=P_2(u^{\perp}),
\]
 i.e.,
\begin{equation}
\label{raz}
P_1(u)+P_0(u)=P_2(u^{\perp}).
\end{equation}
 The relations
  $P_1(u)P_0(u)=0$ and  $P_2(u^{\perp})\subseteq P_0(u),$
imply $P_1(u)P_2(u^{\perp})=0.$  Multiplying both sides of
\eqref{raz} by $P_1(u)$ we obtain
\[
P_1(u)[P_1(u)+P_0(u)]=P_1P_2(u^{\perp}),
\]
 i.e.,
\[
P_1(u)=0.
\]

 (b) Note that if
 $u, v \in \mathcal{GT}$  are orthogonal, then, by virtue of \cite[Lemma 1.8]{fr2}
 we have
 $$
 P_0(u+v)=P_0(u)P_0(v).
 $$
Since $P_1 (u) = P_1 (u^\perp) = \{0\},$  it follows that
\begin{equation}
\label{for}
I = P_2 (u^
\perp) + P_0 (u^\perp)
\end{equation}
and
\[
I =I\cdot I=(P_2(u) + P_0(u))(P_2(u^\perp) + P_0(u^\perp))=\]
\[=P_2 (u)+ P_2(u^\perp),
\]
i.e.,
\begin{equation}
\label{foru}
I = P_2 (u) + P_2 (u^\perp).
\end{equation}
Applying \eqref{for} and \eqref{foru} we obtain $P_2 (u) = P_0
(u^\perp).$ This completes the proof of the lemma.

\begin{lemma}
\label{lem2} If $\nabla$ is a Boolean algebra and $v_i\in\nabla,$
$v_i\diamondsuit v_j$ for $i\neq j, i, j=\overline{1, n},$ then
\[
P_2\left(\sum\limits_{i=1}^{n}v_i\right)=\sum\limits_{i=1}^{n}P_2(v_i).
\]
\end{lemma}

Proof. Since  $v_i\diamondsuit v_j,$ for $i\neq j, i,
j=\overline{1, n},$ it follows by \cite[Corollary 3.4]{fr2} that
for  $i\neq j$
 we have
\[P_2(v_i)P_2(v_j)=0, \]
\[P_2(v_i)P_0(v_j)=P_2(v_i),\]
\[P_0(v_i)P_0(v_j)=P_0(v_i+v_j).\]
By virtue of Lemma~\ref{lem1}~(a) we have  $P_2(v_i)+P_0(v_i)=I.$
Using the relations given above, we obtain
\[I=I^n=\prod\limits_{i=1}^{n}[P_2(v_i)+P_0(v_i)]=
\]
\[
=\sum\limits_{i=1}^{n}P_2(v_i)+ \prod\limits_{i=1}^{n}P_0(v_i)=
\sum\limits_{i=1}^{n}P_2(v_i)+P_0\left(\sum\limits_{i=1}^{n}v_i\right),\]
i.e.,
\[\sum\limits_{i=1}^{n}P_2(v_i)+P_0\left(\sum\limits_{i=1}^{n}v_i\right)=I.\]
The relation
\[P_2\left(\sum\limits_{i=1}^{n}v_i\right)+P_0\left(\sum\limits_{i=1}^{n}v_i\right)=I,\]
implies
 \[
P_2\left(\sum\limits_{i=1}^{n}v_i\right)=\sum\limits_{i=1}^{n}P_2(v_i).
\]
This proves the lemma.

 In what follows, we assume that  $Z$ is a   $PE$
atomic neutral strongly facially symmetric space and there exists
a geometric tripotent $e$  for which $Z_2 (e) =Z.$

Let $\mathrm{T}$  be a maximal family of mutually orthogonal
minimal geometric tripotents from $Z^{\ast},$ i.e.,
$\mathrm{T}=\{v_i: v_i\diamondsuit v_j, i\neq j, i, j\in J\},$
where each $v_i$ is minimal and no minimal geometric tripotent is
orthogonal to all $v_i, i\in J.$ Such a family exists by Zorn's
Lemma.

\begin{lemma}
\label{lem3} If $\mathrm{T}$ is a maximal family of mutually
orthogonal minimal geometric tripotents from $Z^{\ast}$ and
$v=\sup\{v_i: v_i\in J\},$ then $v=e.$
\end{lemma}

Proof.  Suppose that $v\neq e.$ Then  $v^{\perp}\neq 0.$ Since $Z$
 is a $PE$ atomic $SFS$-space,
it follows that there exists a norm exposed point $x$ in
$F_{v^{\perp}}.$ The geometric tripotent $u$ corresponding to  $x$
is minimal. Therefore,  $u\leq v^{\perp}.$ The relation
$v^{\perp}\leq v^{\perp}_i$ implies  $u\leq v^{\perp}_i$ for $i\in
J.$ Hence $u\diamondsuit v_i, i\in J.$ This contradicts the
maximality of the family $\mathrm{T}.$  Thus $v=e,$ which
completes the proof of the lemma.

\begin{lemma}
\label{lem4} If  $\mathrm{T}$ is a maximal family of mutually
orthogonal minimal geometric tripotents from $Z^{\ast},$ then
$\mathrm{T}$ separates the points of  $Z,$ i.e.,  for any $x\in Z,
x\neq 0$ there exists a $v\in \mathrm{T}$ for which $v(x)\neq 0.$
\end{lemma}

Proof. Let  $\mathrm{T}$ be a maximal family of mutually
orthogonal minimal geometric tripotents from $Z^{\ast}.$ According
to Lemma~\ref{lem3} we have $\sup\{v_i: v_i\in J\}=e.$

Suppose that there exists a point  $x\in Z$ with $||x||=1$ such
that
 $v_i(x)=0$ for all   $i\in J.$ The minimality of each  $v_i$ implies the existence of a
 $z_i\in F_{v_i}$ such
that  $P_2(v_i)(x)=v_i(x) z_i$ for any $x\in Z.$ Since $v_i(x)=0,$
 it
follows that  $P_2(v_i)(x)=0,$ and since $P_1(u)=0,$  it follows
that $x\in Z_0(v_i)$ for all $i\in J.$ By virtue of
Lemma~\ref{lem1}~(a)  we have   $P_2(v_i^\perp)=P_0(v_i).$
Therefore, $x\in Z_2(v_i^\perp)$
   for all  $i\in J.$
Take a minimal geometric tripotent $v$ for which   $x\in F_v.$
 We have  $v\diamondsuit v_i$ for all   $i\in
J.$ Thus,  $v\diamondsuit \sup v_i,$ i.e.,
 $v\diamondsuit e,$
 whence $x\in P_0(e).$ It follows from  $P_0(e)=0$ that  $x=0.$
 This contradicts the
assumption $x\neq 0.$  Thus, $\mathrm{T}$ separates the points of
$Z,$ as required.

Let  $\mathrm{T}=\{v_i\}_{i\in J}$
 be a maximal family of mutually orthogonal minimal geometric
tripotents  from $Z^{\ast},$ and let
\[
\ell_1(\mathrm{T})=
\left\{\{\lambda_i\}_{i\in J}: \sum\limits_{i\in J}|\lambda_i|<+\infty\right\}.
\]
Then $\ell_1(\mathrm{T})$  is a Banach space with respect to the
norm
\[
\|\{\lambda_i\}_{i\in J}\|=
\sum\limits_{i\in J}|\lambda_i|.
\]

The following theorem is the main result of this paper; it
describes strongly facially symmetric spaces isometrically
isomorphic to preduals of atomic commutative von Neumann algebras.

\begin{theorem}\label{MTH}
Suppose that $Z$  is a $PE$ neutral atomic strongly facially
symmetric space and there exists a geometric tripotent $e$  for
which  $Z_2 (e) =Z.$ If  $\nabla$   is a Boolean algebra, then $Z$
  is isometrically
isomorphic to the space $\ell_1(\mathrm{T}),$  where
 $\mathrm{T}$ is a maximal family of mutually orthogonal minimal geometric
tripotents from $Z^{\ast}.$
 \end{theorem}

Proof. Let  $z_i\in F_e$ be such that $v_i(z_i)=1.$ Then, for any
$x\in Z$ we have $P_2(v_i)(x)=v_i(x) z_i,$ where $v_i(x)\in
\mathbb{R}, i\in J.$  Let us show that
\[x=\sum\limits_{i\in J}v_i(x)z_i\]
and
\[||x||=\sum\limits_{i\in J}|v_i(x)|. \]
For $i_1,...,i_n\in J,$  we have
\[\sum\limits_{k=1}^{n}|v_{i_{k}}(x)|=
\sum\limits_{k=1}^{n}|v_{i_{k}}(x)|||z_{i_k}||=\]
\[
=[z_{i_k}\diamondsuit z_{i_s}, k\neq s]=
\]
\[
=\left\|\sum\limits_{k=1}^{n}v_{i_{k}}(x)z_{i_k}\right\|=
\left\|\sum\limits_{k=1}^{n} P_2(v_{i_{k}})(x)\right\|=
\]
\[
=\left\|P_2\left(\sum\limits_{k=1}^{n} v_{i_{k}}\right)(x)\right\|\leq ||x||,
\]
i.e.,
\[\sum\limits_{k=1}^{n}|v_{i_{k}}(x)|\leq ||x||.
\]
 Therefore, the series
 $\sum\limits_{i\in J}|v_{i}(x)|$ converges, and
 $\sum\limits_{i\in J}|v_{i}(x)|\leq||x||.$
 Next, for $j\in J$ we have
 \[
 v_j\left(x-\sum  v_{i}(x)z_i\right)=v_j(x)-v_j(x)=0.
 \]
Since  $\mathrm{T}$
 separates the points of  $Z,$ it follows that
 $x=\sum\limits_{i\in J} v_{i}(x)z_i.$
 Thus, the correspondence
 \[
x\in Z\mapsto \{v_{i}(x)\}
 \]
is an isometric isomorphism between  $Z$ and $\ell_1(\mathrm{T}).$
This completes the proof of the theorem.

 \begin{corollary}
Suppose that $Z$  is a real finite-dimensional neutral strongly
facially symmetric space and there exists a geometric tripotent
$e$ for which $Z_2 (e) =Z.$ If  $\nabla$  is a Boolean algebra,
then $Z$ is isometrically isomorphic to the space $\mathbb{R}^{n}$
with norm \[ ||z||=|z_1|+...+|z_n|.
\]
\end{corollary}

\end{document}